\documentclass{article}
\usepackage[utf8]{inputenc}
\usepackage[english]{babel}
\usepackage{amsthm}
\usepackage{graphicx}

\begin{document}

\title{The $n$th+1 Prime Number Limit Formulas}
\author{Artur Kawalec} 

\date{}
\maketitle

\begin{abstract}
A new derivation of Golomb's limit formula for generating the $n$th$+1$ prime number is presented. The limit formula is derived by extracting $p_{n+1}$ from Euler's prime product representation of the Riemann zeta function $\zeta(s)$ in the limit as $s$ approaches infinity. Also, new variations of these limit formulas are explored, such as the logarithm and a half-prime formulas for the $p_{n+1}$ .
\end{abstract}

\section{Introduction}
An interesting set of limit formulas for generating the $n$th$+1$ prime was proved by Golomb in a 1976 paper [1], which since then, did not appear in any reference literature about prime numbers to the best of my knowledge. The essence of these formulas is that $p_{n+1}$ is generated in the limit of the Riemann zeta function and a partial Euler prime product up the $n$th order, thus all primes up to the $nth$ order must be known in order to compute the $nth+1$ prime. Golomb's original proof is based on applying probability distribution theory on natural numbers. In this paper, a simple derivation of Golomb's main formula is presented by inverting the Euler's prime product representation of the Riemann zeta function $\zeta(s)$, and showing that $p_{n+1}$ can be extracted in the limit as $s\to \infty$. We will also derive new variations of these formulas, such as the logarithm formula and the half-prime formula for the $p_{n+1}$ and explore some of their properties.

\section{The Main Formula}
Let $p_n$ be a sequence of prime numbers defined such that $p_1=2$, $p_2=3$, $p_3=5$ and so on, then for $n\geq 0$, the prime number of order $n+1$ is given by the Golomb's limit formula,
\begin{equation}\label{eq:1}
p_{n+1} = \lim_{s\to\infty} \left(1-\frac{Q_n(s)}{\zeta(s)}\right)^{-1/s}
\nonumber
\end{equation}
where the function $Q_n(s) = \prod_{k=1}^n\ \left(1-p^{-s}_k\right)^{-1}$ is defined as the partial Euler prime product up the $n$th order for $n\geq 1$, and $Q_n(s)=1$ for $n=0$.

\begin{proof}

The well known Euler's prime product formula establishes a key connection between the Riemann zeta function and infinite prime number products,

\begin{equation}\label{eq:2}
\zeta(s) = \prod_{n=1}^\infty\ \left(1-\frac{1}{p^s_n}\right)^{-1}
\nonumber
\end{equation}

which is absolutely convergent for $\Re(s)>1$. In order to solve the Euler's prime product formula for $p_1$, one can write,

\begin{equation}\label{eq:3}
\zeta(s) =  \left(1-\frac{1}{p^s_1}\right)^{-1}\epsilon_2(s)
\nonumber
\end{equation}

where $\epsilon_2(s)$ is an infinite prime product starting with $p_2$ term. In general, $\epsilon_k(s)$ is defined as

\begin{equation}\label{eq:4}
\epsilon_k(s) = \prod_{n=k}^\infty\ \left(1-\frac{1}{p^s_n}\right)^{-1}
\nonumber
\end{equation}

By solving for $p_1$ in Equation (3) results in an expression,

\begin{equation}\label{eq:5}
p_1 = \left(1-\frac{\epsilon_2(s)}{\zeta(s)}\right)^{-1/s}
\nonumber
\end{equation}

We then consider taking the limit $s\to\infty$ and note that if the ratio of $\epsilon_2(s)/\zeta(s)$ is expanded into prime products as shown next,

\begin{equation}\label{eq:6}
\frac{\epsilon_2(s)}{\zeta(s)} =\frac{\left(1-p_2^{-s}\right)^{-1}\left(1-p_3^{-s}\right)^{-1}\left(1-p_4^{-s}\right)^{-1}\textnormal{...}}{\left(1-p_1^{-s}\right)^{-1}\left(1-p_2^{-s}\right)^{-1}\left(1-p_3^{-s}\right)^{-1}\textnormal{...}}
\nonumber
\end{equation}

then the higher order products tend to cancel out, which allows the product term due to $p_1$ in the Riemann zeta function to dominate the limit. Furthermore, because $p_{n}^{-s}\gg p_{n+1}^{-s}$, then the higher order prime products tend to converge to $1$ much faster than the leading prime product due to $p_1$ in the Riemann zeta function, as a result, the contribution due to $\epsilon_2(s)$ in the limit becomes negligible, hence $\epsilon_2(s)\to 1$. The limit can be then expressed as,

\begin{equation}\label{eq:7}
p_1 = \lim_{s\to\infty} \left(1-\frac{1}{\zeta(s)}\right)^{-1/s}
\nonumber
\end{equation}

Numerical computation of this limit results in a correct convergence $p_1\to 2$, which is summarized in Table 1 for $s=10$ and $s=100$.
We can then use the same approach to generate the next prime in the sequence, but now we must take into account a known prime product term due to $p_1$. Hence, if we express the Riemann zeta function as,

\begin{equation}\label{eq:8}
\zeta(s) = \left(1-\frac{1}{p^s_1}\right)^{-1}\left(1-\frac{1}{p^s_2}\right)^{-1}\epsilon_3(s)
\nonumber
\end{equation}

Then by solving for $p_2$ and taking the limit $s\to\infty$ yields,

\begin{equation}\label{eq:9}
p_2 = \lim_{s\to\infty} \left(1-\frac{1}{\zeta(s)(1-p_1^{-s})}\right)^{-1/s}
\nonumber
\end{equation}

in which $\epsilon_3(s)\to 1$ for the same reason that higher order primes products tend to cancel out. And so, if we set $p_1=2$ into the above limit, then $p_2 \to3$. We can then further generalize this procedure to compute $p_{n+1}$ by repeating the same pattern. If we define a finite prime product up the $n$th order,

\begin{equation}\label{eq:10}
Q_n(s) = \prod_{k=1}^n\ \left(1-\frac{1}{p^s_k}\right)^{-1}
\nonumber
\end{equation}

for $n\geq 1$ and $Q_n(s) = 1$ for $n=0$, then $p_{n+1}$ can be written as,

\begin{equation}\label{eq:11}
p_{n+1} = \lim_{s\to\infty} \left(1-\frac{Q_n(s)}{\zeta(s)}\right)^{-1/s}
\nonumber
\end{equation}

where $\epsilon_{n+1}(s)\to 1$. This completes the proof.
\end{proof}

Table 1 summarizes computation of the first 10 prime numbers computed using Equation (1) for $s=10$ and $s=100$, while Table 2 summarizes computation of a few higher order primes for $s=1000$. It is readily seen that the limit is converging to the $p_{n+1}$ prime, however it's generally difficult to compute it on a standard floating point arithmetic because of the fast convergence of $Q_n(s)$ and $\zeta(s)$ to $1$. For instance, $\zeta(100)$ can be estimated to be on the order of $1+2^{-100}$ which has 31 trailing zeros after the decimal point. It is therefore more appropriate to use a finite precision arithmetic in computation. In Tables 1 and 2, the $N$ function in Mathematica was used to compute these limits.

Also, in order to get a better insight into this limit, the function $p_{n+1}(s)$ is plotted vs. $s$ in Figure 1 for the first few values of $n$, where we can observe graphically how $p_{n+1}(s)$ converges to $p_{n+1}$ as $s$ increases. The error associated with convergence can understood as the culmination of infinite prime products terms $\epsilon_{n+2}(s)$. This can be shown by writing,

\begin{equation}\label{eq:8}
\frac{Q_n(s)}{\zeta(s)} = \frac{1}{\epsilon_{n+1}(s)}
\end{equation}

and if we express $\epsilon_{n+1}(s)$ as,

\begin{equation}\label{eq:8}
\epsilon_{n+1}(s)=\left(1-\frac{1}{p^{s}_{n+1}}\right)^{-1}  \epsilon_{n+2}(s)
\end{equation}

Then, by substituting this result to Equation (1), one essentially recovers $p_{n+1}$ as

\begin{equation}\label{eq:8}
p_{n+1} = \lim_{s\to\infty}\left(1-\frac{1}{\epsilon_{n+2}(s)}+\frac{1}{p_{n+1}^s\epsilon_{n+2}(s)}\right)^{-1/s}
\end{equation}

 because of the fact that $p_{n}^{-s}\gg p_{n+1}^{-s}$, where the higher order prime products of  $\epsilon_{n+2}(s)$ rapidly approach $1$ as $s \to \infty$ which shows that Equation (1) is convergent to $p_{n+1}$.

\begin{table}[ht]
\caption{Evaluation of $p_{n+1}(s)$ to 15 decimal places by equation 1} 
\centering 
\begin{tabular}{c c c} 
\hline\hline 
n+1 & s=10 & s=100\\ [0.5ex] 
\hline 
1  & 1.996546424130332 & 1.999999999999999 \\ 
2  & 2.998128944738979 & 2.999999999999999 \\
3  & 4.982816482987932 & 4.999999999999999 \\
4  & 6.990872151877531 & 6.999999999999999 \\
5  & 10.795904253794409 & 10.999999993885992 \\
6  & 12.882858209904345 & 12.999999999999709 \\
7  & 16.454690036492369 & 16.999997488242396 \\
8  & 18.700432429563358 & 18.999999999042078 \\
9  & 22.653649208924189 & 22.999999999980263 \\
10 & 27.560268802131417  & 28.999632082761238 \\ [1ex] 
\hline 
\end{tabular}
\label{table:nonlin} 
\end{table}

\begin{table}[ht]
\caption{Evaluation of $p_n(s)$ to 15 decimal places by Equation (1)} 
\centering 
\begin{tabular}{c c c} 
\hline\hline 
n & Expected $p_n$& Calculated $p_n(s)$, s=1000\\ [0.5ex] 
\hline 
$10^1$  & 29 & 28.999999999999999 \\ 
$10^2$  & 541 & 540.99999122731532 \\
$10^3$  & 7919 & 7914.878107364037780 \\
$10^4$  & 104729 & 104488.769372683995648 \\ [1ex] 
\hline 
\end{tabular}
\label{table:nonlin} 
\end{table}

\begin{figure}[h]
  \centering
  \includegraphics[width=130mm]{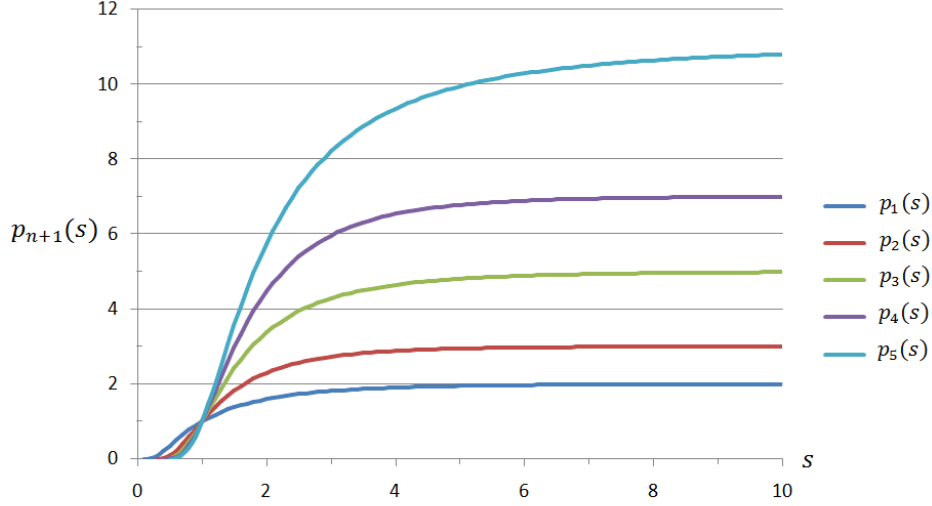}\\
  \caption{Plot of $p_{n+1}(s)$ by Equation (1)}\label{1}
\end{figure}

\section{Additional Formulas}

Many simple variations of Equation (1) can be obtained by algebraic manipulation. For example, if $a$ is real and $a>0$, then $p_{n+1}$ raised to power $a$ is given by the following limit formula,

\begin{equation}\label{eq:12}
p_{n+1}^a = \lim_{s\to\infty} \left(1-\frac{Q_n(as)}{\zeta(as)}\right)^{-1/s}
\end{equation}

\begin{proof}
By substituting $s\to as$ in Equation (1).
\end{proof}

This formula can lead to interesting limit identities. For example, if one sets $n=0$ and $a=1/2$, then one could obtain the square root of $p_1$ as,

\begin{equation}\label{eq:15}
\lim_{s\to\infty} \left(1-\frac{1}{\zeta(s/2)}\right)^{-1/s}=\sqrt2
\end{equation}

Another limit formula which also appears in Golomb's paper with $a=1$ is given by,

\begin{equation}\label{eq:13}
p_{n+1}^a = \lim_{s\to\infty} \left[\zeta(as)-Q_n(as)\right]^{-1/s}
\end{equation}

\begin{proof}

By factoring out $\zeta(s)$ from Equation (1) results in,

\begin{equation}\label{eq:14}
p_{n+1}^a = \lim_{s\to\infty} \zeta(as)^{1/s}\left[\zeta(as)-Q_n(as)\right]^{-1/s}
\end{equation}

It follows that as $s\to\infty$, then $\zeta(as)^{1/s}$ converges to $1$. The resulting limit must then converge to $p_{n+1}^a$.
\end{proof}

Using Equation (17), the square root of $p_1$ can be re-written as,

\begin{equation}\label{eq:15}
\lim_{s\to\infty} \left[\zeta(s/2)-1\right]^{-1/s}=\sqrt2
\end{equation}

The essential mechanism behind these limits is that higher order prime products become negligible as $s\to \infty$, and we are just solving for the dominant prime which is encoded in the Riemann zeta function. In another example, if $n=1$ and $a=3$, then using Equation (17), one could generate a limit identity for the cube of $p_2$,

\begin{equation}\label{eq:16}
\lim_{s\to\infty} \left[\zeta(3s)-(1-1/2^{3s})^{-1}\right]^{-1/s}=27
\end{equation}

Equation (17) can be used to generate a logarithm formula for $p_{n+1}$ as shown next,

\begin{equation}\label{eq:17}
\log(p_{n+1})=\lim_{s\to\infty} -\frac{Q_{n}^{'}(s)-\zeta^{'}(s)}{Q_{n}(s)-\zeta(s)}
\end{equation}

\begin{proof}
By applying a logarithm to Equation (17), and by differentiating both sides with respect to $a$ and absorbing $as\to s$ in the limit completes the proof.
\end{proof}

As an example, let us consider the simplest case for $n=0$. Then the partial Euler products results in $Q_0(s)=1$ and $Q^{'}_0(s)=0$, hence we get another limit identity,

\begin{equation}\label{eq:19}
\lim_{s\to\infty} \frac{\zeta^{'}(s)}{1-\zeta(s)}=\log(2)
\end{equation}

in which the logarithm of the first prime is expressed in the limit of the Riemann zeta function and its first derivative. A numerical computation of this limit with $s=100$ results in convergence to $0.693147180$.

Another limit formula variant is given by,

\begin{equation}\label{eq:20}
p_{n+1}^{a}=\lim_{s\to\infty} \left(\frac{Q_n^{'}(as)}{Q_n(as)}-\frac{\zeta^{'}(as)}{\zeta(as)}\right)^{-1/s}
\end{equation}

\begin{proof}
This formula is determined by differentiation of the logarithm of Equation (15) with respect to $a$,

\begin{equation}\label{eq:22}
\log(p_{n+1})=\lim_{s\to\infty} \left(\frac{Q_n(as)}{\zeta(as)}\right)\left(1-\frac{Q_n(as)}{\zeta(as)}\right)^{-1}\left(\frac{Q_n^{'}(as)}{Q_n(as)}-\frac{\zeta^{'}(as)}{\zeta(as)}\right)
\end{equation}

which can be further simplified by noticing that the first term approaches $1$, and the second term approaches $p_{n+1}^{as}$ due to Equation (15). And so, by re-arranging the terms we get,

\begin{equation}\label{eq:23}
p_{n+1}^{a}=\lim_{s\to\infty} \left(\frac{1}{\log(p_{n+1})}\right)^{-1/s}\left(\frac{Q_n^{'}(as)}{Q_n(as)}-\frac{\zeta^{'}(as)}{\zeta(as)}\right)^{-1/s}
\end{equation}

And since the logarithm term tends to approach $1$, this complete the proof.
\end{proof}

As a simple result of this formula, if we let $n=0$, then one can recover an asymptotic relation for,

\begin{equation}\label{eq:23}
-\frac{\zeta'(s)}{\zeta(s)}\sim 2^{-s}
\end{equation}

This formula is also interesting because of its relation to $\zeta'(s)/\zeta(s)$ term which often appears in the study of prime numbers. For instance,
\begin{equation}\label{eq:23}
-\frac{\zeta'(s)}{\zeta(s)} = \sum_{n=1}^{\infty} \frac{\Lambda(n)}{n^s}
\end{equation}

where $\Lambda(n)$ is the von Mangoldt's function.

\begin{figure}[h]
  \centering
  \includegraphics[width=130mm]{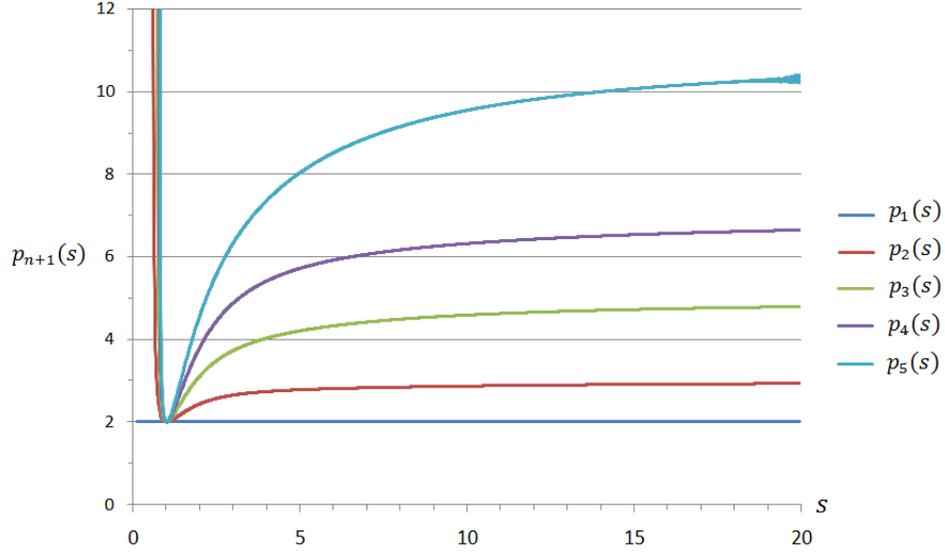}\\
  \caption{Plot of $p_{n+1}(s)$ by Equation (28)}\label{2}
\end{figure}

We can also derive a half-prime limit formula based only first derivatives of $Q_n(s)$ and $\zeta(s)$,

\begin{equation}\label{eq:25}
\frac{p_{n+1}}{2} = \lim_{s\to\infty}\left(1-\frac{Q^{'}_n(s)}{\zeta^{'}(s)}\right)^{-1/s}
\end{equation}

where $n\geq 0$.

\begin{proof}
This formula can be determined by factoring out $\zeta^{'}(s)/\zeta(s)$ from Equation (21) which results in,
\begin{equation}\label{eq:26}
\log(p_{n+1})=\lim_{s\to\infty}-\frac{\zeta^{'}(s)}{\zeta(s)}\left(1-\frac{Q_{n}(s)}{\zeta(s)}\right)^{-1}\left(1-\frac{Q_{n}^{'}(s)}{\zeta^{'}(s)}\right)
\end{equation}
Using Equation (15) and noting that the asymptotic relation for $-\zeta^{'}(s)/\zeta(s)\to 2^{-s}$, we can write the above result as,

\begin{equation}\label{eq:27}
\log(p_{n+1})=\lim_{s\to\infty} 2^{-s}p_{n+1}^{-s}\left(1-\frac{Q_{n}^{'}(s)}{\zeta^{'}(s)}\right)
\end{equation}

And by re-arranging the terms we get,

\begin{equation}\label{eq:28}
p_{n+1}=\lim_{s\to\infty} \left(\frac{2^{-s}}{\log(p_{n+1})}\right)^{-1/s}\left(1-\frac{Q_{n}^{'}(s)}{\zeta^{'}(s)}\right)^{-1/s}
\end{equation}

In the limit as $s\to \infty$, the first term approaches $2$, and so this leads to,

\begin{equation}\label{eq:29}
\frac{p_{n+1}}{2} = \lim_{s\to\infty} \left(1-\frac{Q^{'}_n(s)}{\zeta^{'}(s)}\right)^{-1/s}
\end{equation}

\end{proof}

This formula tends to converge much slower to $p_{n+1}$ than the limit formula of Equation (1). For example, if $n=3$ and $s=100$, then $p_{n+1}$ evaluated by Equation (1) gives $6.999999999$, while $p_{n+1}$ evaluated by Equation (28) gives $6.928114662$. Also, the function $p_{n+1}(s)$ of Equation (28) is plotted vs.$s$ in Figure 2 for the first few values of $n$.  We note that for $n\geq 1$, the function $p_{n+1}(s)$ sharply drops from infinity approaching a value of $2$, and then slowly settles to $p_{n+1}$ as $s$ increases.

\section{Conclusion}
We have shown a simple proof of Golomb's main formula without rigorous analysis and have validated it numerically for many cases. We have also found additional variants of these formulas such as the logarithm and half-prime formulas for the next prime. Golomb in [1] mentioned that more formulas could be found.

\section{Acknowledgement}
I would like to thank Igor Acimovic for reviewing the paper and providing valuable suggestions.

\end{document}